\documentclass[12pt]{article}
\usepackage{latexsym,amssymb,amsmath,amsfonts,pictex,enumerate,amsthm,dsfont}
\usepackage{tikz}
\usetikzlibrary{decorations.markings}
\usepackage{scalerel,stackengine}
\usepackage[active]{srcltx}

\theoremstyle{definition}
\newtheorem{theorem}{Theorem}

\newtheorem{lemma}[theorem]{Lemma}

\stackMath
\newcommand\reallywidehat[1]{%
\savestack{\tmpbox}{\stretchto{%
  \scaleto{%
    \scalerel*[\widthof{\ensuremath{#1}}]{\kern-.6pt\bigwedge\kern-.6pt}%
    {\rule[-\textheight/2]{1ex}{\textheight}}
  }{\textheight}%
}{0.5ex}}%
\stackon[1pt]{#1}{\tmpbox}%
}
\begin{document}

\title{Non-tangential limits for analytic Lipschitz functions
}


\author{Stephen Deterding    
\thanks{
              University of Kentucky, Lexington, KY 40506, USA,  Email:stephen.deterding@westliberty.edu, Present address: West Liberty University, West Liberty, WV 26074, USA
            }    
}
\date{}



\maketitle

\begin{abstract}
Let $U$ be a bounded open subset of the complex plane. Let $0<\alpha<1$ and let $A_{\alpha}(U)$ denote the space of functions that satisfy a Lipschitz condition with exponent $\alpha$ on the complex plane, are analytic on $U$ and are such that for each $\epsilon >0$, there exists $\delta >0$ such that for all $z$, $w \in U$, $|f(z)-f(w)| \leq \epsilon |z-w|^{\alpha}$ whenever $|z-w| < \delta$. We show that if a boundary point $x_0$ for $U$ admits a bounded point derivation for $A_{\alpha}(U)$ and $U$ has an interior cone at $x_0$ then one can evaluate the bounded point derivation by taking a limit of a difference quotient over a non-tangential ray to $x_0$. Notably our proofs are constructive in the sense that they make explicit use of the Cauchy integral formula.
\end{abstract}

\section{Background and statement of results}
In this paper, we consider the behavior of Lipschitz functions which are analytic on a bounded open subset of the complex plane and how much analyticity extends to the boundary of the domain. Let $U$ be an open subset in the complex plane and let $0 < \alpha < 1$. A function $f : U \to \mathbb{C}$ satisfies a Lipschitz condition with exponent $\alpha$ on $U$ if there exists $k>0$ such that for all $z,w \in U$
\begin{equation}
\label{lip condition}
|f(z) - f(w)| \leq k |z-w|^{\alpha}
\end{equation}

\bigskip

Let Lip$\alpha(U)$ denote the space of functions that satisfy a Lipschitz condition with exponent $\alpha$ on $U$. Lip$\alpha(U)$ is a Banach space with norm given by $||f||_{Lip\alpha(U)} = \sup_{U} |f| + k(f)$, where $k(f)$ is the smallest constant that satisfies \eqref{lip condition}. If we let $||f||_{Lip\alpha(U)}' = k(f)$ then $||f||_{Lip\alpha(U)}'$ is a seminorm on Lip$\alpha(U)$. 

\bigskip

An important subspace of Lip$\alpha(U)$ is the little Lipschitz class, lip$\alpha(U)$, which consists of those functions in Lip$\alpha(U)$ that also satisfy the additional property that for each $\epsilon >0$, there exists $\delta >0$ such that for all $z$, $w$ in $U$, $|f(z)-f(w)| \leq \epsilon |z-w|^{\alpha}$ whenever $|z-w| < \delta$. 

\bigskip

The importance of lip$\alpha(U)$ is illustrated by the following result of De Leeuw \cite{DeLeuuw}. Let $\Delta$ be a closed disk. Then the restriction spaces Lip$\alpha(\Delta) = \{f|\Delta: f \in \textnormal{ Lip} \alpha(\mathbb{C})\}$ and lip$\alpha(\Delta) = \{f|\Delta: f \in \textnormal{ lip} \alpha(\mathbb{C})\}$ are Banach spaces and lip$\alpha^{**}(\Delta)$ is isometrically isomorphic to Lip$\alpha(\Delta)$. Thus the weak-star topology can be applied to Lip$\alpha(\Delta)$ as the dual of lip$\alpha^{*}(\Delta)$.

\bigskip

Let $U$ be a bounded open subset of the complex plane. We will restrict our study to those functions in lip$\alpha(\mathbb{C})$ which are analytic on $U$. Let $A_{\alpha}(U) = \{f \in \textnormal{lip}\alpha: \overline{\partial} f = 0 \textnormal{ on } U\}$, where $\overline{\partial}f = \dfrac{1}{2} \left(\dfrac{\partial f}{\partial x} + i \dfrac{\partial f}{\partial y} \right)$. For an arbitrary set $E \subset \mathbb{C}$, let $A_{\alpha}(E) = \bigcup \{A_{\alpha}(U): U \textnormal{ open }, E \subset U\}$.

\bigskip

While the functions in $A_{\alpha}(U)$ are differentiable on the interior of $U$, they need not be differentiable on the boundary of $U$. In this paper, we consider the question of how close the functions in $A_{\alpha}(U)$ come to being differentiable at boundary points of $U$. To answer this question we will make use of the concept of a bounded point derivation. For $x_0 \in \mathbb{C}$, it is known that $A_{\alpha}(U \cup \{x_0\})$ is dense in $A_{\alpha}(U)$. \cite[Lemma 1.1]{Lord} Thus we say that $A_{\alpha}(U)$ admits a bounded point derivation at $x_0$ if the map $f\to f'(x_0)$ extends from $A_{\alpha}(U \cup \{x_0\})$ to a bounded linear functional on $A_{\alpha}(U)$. Equivalently, $A_{\alpha}(U)$ admits a bounded point derivation at $x_0$ if and only if there exists a constant $C >0$ such that 

\begin{equation}
\label{a(U)bpd}
|f'(x_0)| \leq C ||f||_{Lip\alpha(\mathbb{C})},
\end{equation}

\bigskip

\noindent for all $f$ in $A_{\alpha}(U \cup \{x_0\})$.

\bigskip

The existence of a bounded point derivation at $x_0$ shows that the functions in $A_{\alpha}(U)$ possess some semblance of analytic structure at $x_0$. If, in addition, $U$ has an interior cone at $x_0$, a more explicit description of this analytic structure can be obtained. We say that $U$ has an interior cone at $x_0$ if there is a segment $J$ ending at $x_0$ and a constant $k > 0$ such that dist$(x, \partial U) \geq k |x-x_0|$ for all $x$ in $J$. The segment $J$ is called a non-tangential ray to $x_0$. It is a result of O'Farrell \cite{O'Farrell2014} that if $U$ has an interior cone at a boundary point $x_0$, then a bounded point derivation on $A_{\alpha}(U)$ at $x_0$ can be evaluated by taking the limit of the difference quotient over a non-tangential ray to $x_0$.   To be precise, O'Farrell has proven the following theorem.

\begin{theorem}
\label{lip}

Let $0 < \alpha <1$, and let $U$ be an open set with $x_0$ in $\partial U$. Suppose that $U$ has an interior cone at $x_0$ and that $J$ is a non-tangential ray to $x_0$. If $A_{\alpha}(U)$ admits a bounded point derivation $D$ at $x_0$, then for every $f$ in $A_{\alpha}(U)$,

\begin{equation*}
\label{lipeq}
Df = \lim_{x \to x_0, x \in J} \dfrac{f(x)-f(x_0)}{x-x_0}.  
\end{equation*}

\end{theorem}

\bigskip

Thus the difference quotient for boundary points that admit bounded point derivations for $A_{\alpha}(U)$ converges when taken over a non-tangential ray to the point. This illustrates the additional analytic structure of functions in $A_{\alpha}(U)$ at these points.

\bigskip

 O'Farrell comments that the methods used in his proof of Theorem \ref{lip} are nonconstructive, involving abstract measures and duality arguments from functional analysis as opposed to using the Cauchy integral formula directly, and suggests that it should be possible to give a proof using constructive techniques. In this paper we present a constructive proof of Theorem \ref{lip}, which confirms O'Farrell's conjecture. In Section $2$ we review some key properties of $A_{\alpha}(U)$ and in Section $3$ we prove Theorem \ref{lip} using constructive techniques.

\section{Preliminary Results}

We begin by reviewing the Hausdorff content of a set, which is defined using measure functions. A measure function is a monotone nondecreasing function $h : [0, \infty) \to [0, \infty)$. For example, $r^{\beta}$ is a measure function for $0 \leq  \beta < \infty$. If $h$ is a measure function then the Hausdorff content $M_h$ associated to $h$ is defined by 

\[ M_h(E) = \inf \sum h(\textnormal{diam } B), \]

\bigskip

\noindent where the infimum is taken over all countable coverings of $E$ by balls and the sum is taken over all the balls in the covering. If $h(r) = r^{\beta}$ then we denote $M_h$ by $M^{\beta}$. The lower $1+ \alpha$ dimensional Hausdorff content $M^{1+\alpha}_*(E)$ is defined by 

\[ M^{1+\alpha}_*(E) = \sup M_h(E), \]

\bigskip

\noindent where the supremum is taken over all measurable functions $h$ such that $h(r) \leq r^{1+\alpha}$ and $r^{-1-\alpha}h(r)$ converges to $0$ as $r$ tends to $0$. The lower $1+ \alpha$ dimensional Hausdorff content is a monotone set function; i.e. if $E \subseteq F$ then $M^{1+\alpha}_*(E) \leq M^{1+\alpha}_*(F)$.

\bigskip

In \cite{Lord}, Lord and O'Farrell gave necessary and sufficient conditions for the existence of bounded point derivations on $A_{\alpha}(U)$ in terms of Hausdorff contents. There are similar conditions for bounded point derivations defined on other function spaces. (\cite{Hallstrom}, \cite{Hedberg})

\begin{theorem}
\label{Lord}
Let $U$ be an open subset of the complex plane with $x_0$ on the boundary of $U$. Let $0<\alpha < 1$. Then $A_{\alpha}(U)$ has a bounded point derivation at $x_0$ if and only if

\begin{equation*}
\sum_{n=1}^{\infty} 4^n M_*^{1+\alpha}(A_n(x_0) \setminus U) < \infty.
\end{equation*}

\end{theorem}

\bigskip


Another key lemma is the following Cauchy theorem for Lipschitz functions which also appears in the paper of Lord and O'Farrell \cite[pg.110]{Lord}.

\begin{lemma}
\label{Cauchy2}

Let $\Gamma$ be a piecewise analytic curve bounding a region $\Omega \in \mathbb{C}$, and suppose that $\Gamma$ is free of outward pointing cusps. Let $0 < \alpha < 1$ and suppose that $f \in$ lip$\alpha(\mathbb{C})$. Then there exists a constant $\kappa >0$ such that

\begin{equation*}
\left| \int f(z)dz\right| \leq \kappa \cdot M_{*}^{1+\alpha}(\Omega \cap S) \cdot ||f||'_{Lip\alpha(\Omega)}.
\end{equation*}

\bigskip

\noindent The constant $\kappa$ only depends on $\alpha$ and the equivalence class of $\Gamma$ under the action of the conformal group of $\mathbb{C}$. In particular this means that $\kappa$ is the same for any curve obtained from $\Gamma$ by rotation or scaling.

\end{lemma}






\section{The proof of the main theorem}


To prove Theorem \ref{lip}, we first note that by translation invariance we may suppose that $x_0 = 0$. Moreover by replacing $f$ by $f-f(0)$ if needed, we may suppose that $f(0) = 0$. In addition, we may suppose that $U$ is contained in the unit disk. Let $J$ be a non-tangential ray to $x_0$ and for each $x$ in $J$, define a linear functional $L_x$ by $L_x(f) = \dfrac{f(x)}{x} - Df$. Then to prove Theorem \ref{lip} it suffices to show that $L_x$ tends to the $0$ functional as $x \to 0$ through $J$. We make the following claim.

\begin{lemma}
\label{claim}
 The collection $\{L_x : x\in J\}$ is a family of bounded linear functionals on $A_{\alpha}(U)$; that is there exists a constant $C >0$ that does not depend on $x$ or $f$ such that $|L_x(f)| \leq C ||f||_{Lip \alpha(\mathbb{C})}$ for all $f$ in $A_{\alpha}(U)$ and all $x \in J$.

\end{lemma}

\bigskip






\bigskip

\begin{proof}

We will first prove Lemma \ref{claim} for the case when $f$ belongs to $A_{\alpha}(U \cup \{0\})$ and then extend to the general case. It follows from \eqref{a(U)bpd} that it is enough to show that $\left|\dfrac{f(x)}{x}\right| \leq C ||f||_{Lip\alpha(\mathbb{C})}$ where the constant $C$ does not depend on $f$ or $x$. If $f$ belongs to $A(U \cup \{0\})$, then there is a neighborhood $\Omega$ of $0$ such that $f$ is analytic on $\Omega$. We can further suppose that $U \subseteq \Omega$. Let $B_n$ denote the ball centered at $0$ with radius $2^{-n}$. Then there exists an integer $N>0$ such that $\Omega$ contains $B_N$ and hence $f$ is analytic inside the ball $B_N$. In addition, there exists an integer $M$ such that $\Omega \subseteq B_M$. Since $J$ is a non-tangential ray to $x_0$, it follows that there is a sector in $\mathring{U}$ with vertex at $x_0$ that contains $J$. Let $C$ denote this sector. It follows from the Cauchy integral formula that

\begin{equation*}
\dfrac{f(x)}{x} = \frac{1}{2 \pi i} \int_{\partial(C \bigcup B_N)} \dfrac{f(z)}{z(z-x)}dz  
\end{equation*}

\bigskip

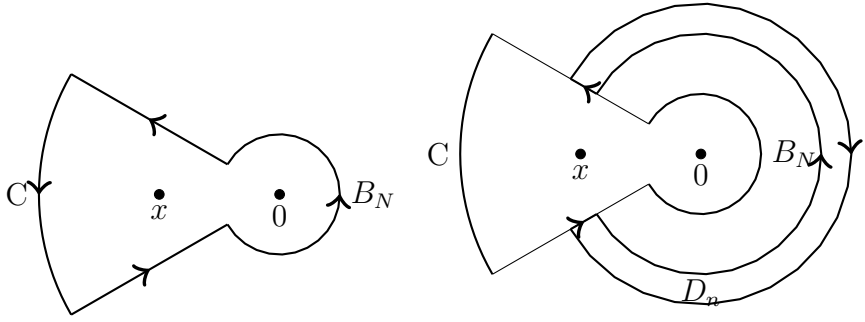
\begin{figure}
\begin{tikzpicture}[scale=0.8]

\tikzset{->-/.style={decoration={
  markings,
  mark=at position .5 with {\arrow[scale=2]{>}}},postaction={decorate}}}

\draw [->-,thick,domain=-150:150] plot ({cos(\x)}, {sin(\x)});
   \draw [->-, thick,domain=150:210] plot ({4* cos(\x)}, { 4*sin(\x)});
   
   \draw [->-, thick](-0.866025, 0.5) to (-3.464102, 2);
   \draw [->-, thick] (-3.464102, -2) to (-0.866025, -0.5);
   \filldraw[fill=black, draw=black] (0,0) circle (.08 cm);
   \node [below] (x) at (0, 0){$0$};
   \node [left] (x) at (-4,0){C};
   \node [right] (x) at (1,0){$B_N$};
   \filldraw[fill=black, draw=black] (-2,0) circle (.08 cm);
   \node [below] (x) at (-2, 0){$x$};





\end{tikzpicture}
\begin{tikzpicture}[scale=0.8]

\tikzset{->-/.style={decoration={
  markings,
  mark=at position .5 with {\arrow[scale=2]{>}}},postaction={decorate}}}

\draw [thick,domain=-150:150] plot ({cos(\x)}, {sin(\x)});
   \draw [  thick,domain=150:210] plot ({4* cos(\x)}, { 4*sin(\x)});
    \draw [ ->-, thick,domain=-150:150] plot ({2* cos(\x)}, { 2*sin(\x)});
     \draw [ ->-, thick,domain=-150:150] plot ({2.5* cos(-\x)}, { 2.5*sin(-\x)});
     
     \draw [->-, thick](-1.732051, 1) to (-2.1650635, 1.25);
     
     \draw [->-, thick] (-2.1650635, -1.25) to (-1.732051, -1);
   \draw (-0.866025, 0.5) to (-3.464102, 2);
   \draw  (-3.464102, -2) to (-0.866025, -0.5);
   \filldraw[fill=black, draw=black] (0,0) circle (.08 cm);
   \node [below] (x) at (0, 0){$0$};
   \node [below] (x) at (0, -1.9){$D_n$};
   \node [left] (x) at (-4,0){C};
   \node [right] (x) at (1,0){$B_N$};
   \filldraw[fill=black, draw=black] (-2,0) circle (.08 cm);
   \node [below] (x) at (-2, 0){$x$};





\end{tikzpicture}

\caption{The contour of integration}
\label{fig3}
\end{figure}

\noindent where the boundary is oriented so that the interior of $C \bigcup B_N$ lies always to the left of the path of integration. (See Figure \ref{fig3}.) Let $D_n = A_n \setminus C$. Then 

\[\dfrac{f(x)}{x} = \frac{1}{2 \pi i} \sum_{n=M}^N \int_{\partial D_n} \dfrac{f(z)}{z(z-x)}dz + \frac{1}{2 \pi i} \int_{|z|=2^{-M}} \dfrac{f(z)}{z(z-x)}dz. \]

\bigskip




\noindent Since $x$ lies on $J$, which is a non-tangential ray to $x_0$, there exists a constant $k>0$ such that for $z \notin U$, $\dfrac{|x|}{|z-x|} \leq k^{-1}$. Thus for $z \notin U$, $\dfrac{|z|}{|z-x|} \leq 1 + \dfrac{|x|}{|z-x|} \leq 1 + k^{-1}$. Hence $\dfrac{1}{|z|\cdot |z-x|} \leq \dfrac{1+ k^{-1}}{|z|^2}$ and therefore

\begin{equation}
\label{diffquot2}
\dfrac{|f(x)|}{|x|} \leq \frac{1}{2 \pi } \sum_{n=M}^N \left|\int_{\partial D_n} \dfrac{f(z)}{z(z-x)}dz\right|+ \dfrac{4^M(1+k^{-1})}{2 \pi} ||f||_{\infty}.
\end{equation}

\bigskip

\noindent Since $\dfrac{f(z)}{z(z-x)}$ is analytic on $D_n \setminus U$ for $M \leq n \leq N$, an application of Lemma \ref{Cauchy2} shows that 

\begin{equation}
\label{intbound2}
\left| \int_{\partial D_n} \dfrac{f(z)}{z(z-x)}dz \right| \leq \kappa M_*^{1+\alpha}(D_n \setminus U) \cdot \left \Vert \dfrac{f(z)}{z(z-x)} \right \Vert _{Lip\alpha(D_n)}'.
\end{equation}

\bigskip

\noindent Recall that the constant $\kappa$ is the same for curves in the same equivalence class. Since the regions $D_n$ differ from each other by a scaling it follows that $\kappa$ doesn't depend on $n$ in \eqref{intbound2}. 

\bigskip

We now show that $\left \Vert \dfrac{f(z)}{z(z-x)} \right \Vert _{Lip\alpha(D_n)}'$ can be bounded by a constant independent of $f$ and $x$. It follows from the definition of the Lipschitz seminorm that

\begin{align*}
\begin{split}
\left \Vert \dfrac{f(z)}{z(z-x)} \right \Vert _{Lip\alpha(D_n)}' &= \sup_{z \neq w; z,w \in D_n} \dfrac{\left|\dfrac{f(z)}{z(z-x)} - \dfrac{f(w)}{w(w-x)}\right|}{|z-w|^{\alpha}}\\[.2in]
&= \sup_{z \neq w; z,w \in D_n} \dfrac{|w(w-x)f(z)-z(z-x)f(w)|}{|z|\cdot |z-x|\cdot|w|\cdot|w-x|\cdot|z-w|^{\alpha}}.
\end{split}
\end{align*}

\bigskip

\noindent Thus it follows from the triangle inequality that

\begin{align}
\label{equationx}
\begin{split}
\left \Vert \dfrac{f(z)}{z(z-x)} \right \Vert _{Lip\alpha(D_n)}' &\leq \sup_{z \neq w; z,w \in D_n} \dfrac{|w(w-x)f(z)-w(w-x)f(w)|}{|z|\cdot |z-x|\cdot|w|\cdot|w-x|\cdot|z-w|^{\alpha}} \\[.2in]
&+ \sup_{z \neq w; z,w \in D_n} \dfrac{|w(w-x)f(w) -z(z-x)f(w)|}{|z|\cdot |z-x|\cdot|w|\cdot|w-x|\cdot|z-w|^{\alpha}}.
\end{split}
\end{align}

\bigskip

\noindent We first bound the first term on the right of \eqref{equationx}

\begin{align*}
&\sup_{z \neq w; z,w \in D_n} \dfrac{|w(w-x)f(z)-w(w-x)f(w)|}{|z|\cdot |z-x|\cdot|w|\cdot|w-x|\cdot|z-w|^{\alpha}} \\[.2in]
&\leq  \sup_{z \in D_n} \dfrac{1}{|z|\cdot |z-x|} \cdot ||f||'_{Lip\alpha(D_n)}.
\end{align*}

\bigskip

\noindent Since $z \notin U$, $\dfrac{1}{|z|\cdot |z-x|} < \dfrac{1+k^{-1}}{|z|^2}$, and therefore,
 
 \begin{align}
 \label{equationy}
 \begin{split}
\sup_{z \neq w; z,w \in D_n} \dfrac{|w(w-x)f(z)-w(w-x)f(w)|}{|z|\cdot |z-x|\cdot|w|\cdot|w-x|\cdot|z-w|^{\alpha}} \leq C 4^n ||f||'_{Lip\alpha(D_n)}.
\end{split}
\end{align}

\bigskip

\noindent We now bound the second term on the right side of \eqref{equationx}. Since $f(0) =0$ it follows that for $w \in \mathbb{C}$, $ \dfrac{|f(w)|}{|w|^{\alpha}} \leq ||f||'_{Lip\alpha(\mathbb{C})}$. Moreover, a computation shows that $w(w-x)-z(z-x) = (w-z)(z+w-x)$. Hence 




\begin{align}
\label{mediumbound}
\begin{split}
&\sup_{z \neq w; z,w \in D_n} \dfrac{|w(w-x)f(w) -z(z-x)f(w)|}{|z|\cdot |z-x|\cdot|w|\cdot|w-x|\cdot|z-w|^{\alpha}} \leq \\[.2in]
 &\left(  \sup_{z \neq w; z,w \in D_n}  \dfrac{|w-z|^{1-\alpha} }{ |z-x|\cdot|w|^{1-\alpha}\cdot|w-x|}  +  \dfrac{|w-z|^{1-\alpha} }{|z|\cdot |z-x|\cdot|w|^{1-\alpha} }\right)   \cdot ||f||'_{Lip\alpha(\mathbb{C})}.
\end{split}
\end{align}

\bigskip

\noindent  Since $x$ lies on $J$, there exists a constant $k>0$ such that $\dfrac{1}{|z-x|} < \dfrac{1+k^{-1}}{|z|}$ and $\dfrac{1}{|w-x|} < \dfrac{1+k^{-1}}{|w|}$. Hence 

\begin{align}
\label{est1}
\begin{split}
\sup_{z \neq w; z,w \in D_n} \dfrac{|w-z|^{1-\alpha} }{ |z-x|\cdot|w|^{1-\alpha}\cdot|w-x|} \leq C \dfrac{2^n \cdot (2^n)^{2-\alpha} }{(2^n)^{1-\alpha}} = C 4^n,
\end{split}
\end{align}

\bigskip

\noindent and 

\begin{align}
\label{est2}
\begin{split}
\sup_{z \neq w; z,w \in D_n} \dfrac{|w-z|^{1-\alpha} }{|z|\cdot |z-x|\cdot|w|^{1-\alpha} } \leq C \dfrac{4^n \cdot (2^n)^{1-\alpha} }{(2^n)^{1-\alpha}} = C 4^n.
\end{split}
\end{align}

\bigskip

\noindent  Then \eqref{mediumbound}, \eqref{est1}, and \eqref{est2} yield 

\begin{equation}
\label{equationz}
\sup_{z \neq w; z,w \in D_n} \dfrac{|w(w-x)f(w) -z(z-x)f(w)|}{|z|\cdot |z-x|\cdot|w|\cdot|w-x|\cdot|z-w|^{\alpha}} \leq C 4^n ||f||'_{Lip\alpha(\mathbb{C})},
\end{equation}

\bigskip

\noindent and it follows from \eqref{equationx}, \eqref{equationy}, and \eqref{equationz} that

\begin{equation}
\label{equationw}
\left \Vert \dfrac{f(z)}{z(z-x)} \right \Vert _{Lip\alpha(D_n)}' \leq C 4^n ||f||'_{Lip\alpha(\mathbb{C})}.
\end{equation}

\bigskip

\noindent Thus \eqref{diffquot2}, \eqref{intbound2}, and \eqref{equationw} together yield

\begin{equation*}
\dfrac{|f(x)|}{|x|} \leq C \sum_{n=1}^{\infty} 4^n M_*^{1+\alpha}(D_n \setminus U) \cdot||f||'_{Lip\alpha(\mathbb{C})}.
\end{equation*}

\bigskip

\noindent Since Hausdorff content is monotone, $M_*^{1+\alpha}(D_n \setminus U) \leq M_*^{1+\alpha}(A_n \setminus U)$ and hence

\begin{equation*}
\dfrac{|f(x)|}{|x|} \leq C \sum_{n=1}^{\infty} 4^n M_*^{1+\alpha}(A_n \setminus U) \cdot||f||_{Lip\alpha(\mathbb{C})},
\end{equation*}

\bigskip

\noindent and it follows from Theorem \ref{Lord} that

\begin{equation*}
\dfrac{|f(x)|}{|x|} \leq C ||f||_{Lip\alpha(U)},
\end{equation*}

\bigskip

\noindent where $C$ does not depend on $x$ or $f$. Thus $L_x(f) \leq C ||f||_{Lip\alpha(\mathbb{C})}$ for $f \in A_{\alpha}(U \cup \{0\})$ and since $A_{\alpha}(U \cup \{0\})$ is dense in $A_{\alpha}(U)$, it follows that $L_x$ is a family of uniformly bounded linear functionals on $A_{\alpha}(U)$. 


\qed
\end{proof}

\bigskip

 To complete the proof of Theorem \ref{lip}, since $A_{\alpha}(U \cup {0})$ is dense in $A_{\alpha}(U)$, there exists a sequence $\{f_j\}$ in $A_{\alpha}(U \cup {0})$ such that $f_j \to f$ in the Lipschitz norm. Since each $f_j$ is analytic in a neighborhood of $0$ and since $Df_j = f_j'(0)$, it follows that for each $j$, $L_x(f_j) \to 0$ as $x \to 0$. It follows from the claim that $|L_x(f)-L_x(f_j)| \leq C ||f-f_j||_{Lip \alpha(U)}$. By first choosing $j$ sufficiently large, the right hand side can be made arbitrarily small. Then by choosing $x$ sufficiently close to $0$, $L_x(f_j)$ can be made arbitrarily close to $0$. Thus $L_x(f) \to 0$ as $x\to 0$ through $J$, which proves Theorem \ref{lip}.




\end{document}